\documentclass[a4paper,11pt]{amsart}

\usepackage{amsmath}
\usepackage{amssymb}
\usepackage{amstext}
\usepackage{amsbsy}
\usepackage{amsopn}
\usepackage{amsthm}
\usepackage{amscd}
\usepackage{amsxtra}
\usepackage{amsfonts}

\usepackage{ifthen}
\usepackage{xspace}
\usepackage{fancyheadings}
\usepackage{layout}
\usepackage{hyperref}

\input xy
\xyoption{matrix}
\xyoption{arrow}
\xyoption{curve}
\newdir{ (}{{}*!/-5pt/@^{(}}
\newcommand{\xycenter}[1]{\begin{center}
                          \mbox{\xymatrix{#1}}
                          \end{center}
                         }

\newboolean{xlabels}
\newcommand{\xlabel}[1]{
                        \label{#1}
                        \ifthenelse{\boolean{xlabels}}
                                   {\marginpar{#1}}
                                   {}
                       }

\newcommand{\ZZ}{\mathbb{Z}}

\newcommand{\pz}{\mathbb{P}}

\newcommand{\FZ}{\mathbb{F}}
\newcommand{\AZ}{\mathbb{A}}

\newcommand{\PP}{\mathbb{P}}

\newcommand{\FF}{\mathbb{F}}

\newcommand{\sB}{{\mathcal B}}

\newcommand{\sN}{{\mathcal N}}
\newcommand{\sO}{{\mathcal O}}

\newcommand{\suchthat}{\, | \,}

\newboolean{probleme}
\newcommand{\problem}[1]
           {\ifthenelse{\boolean{probleme}}
                       {{\bf(PROBLEM: #1)\bf}}
                       {}
           }
\newboolean{zukuenftiges}
\newcommand{\zukunft}[1]
           {\ifthenelse{\boolean{zukuenftiges}}
                       {{\bf(AUSBAUM\"OGLICHKEIT: #1)\bf}}
                       {}
           }
\newboolean{extras}
\newcommand{\extra}[1]
           {\ifthenelse{\boolean{extras}}
                       {{\bf EXTRA #1 EXTRA\bf}}
                       {}
           }

\newboolean{ignore}
\newcommand{\ignore}[1]
           {\ifthenelse{\boolean{ignore}}
                       {{\bf IGNORE #1 IGNORE\bf}}
                       {}
           }

\DeclareMathOperator{\rank}{rank}

\DeclareMathOperator{\map}{map}

\theoremstyle{plain}
\begingroup

\newtheorem{thm}{Theorem}
\newtheorem{cor}[thm]{Corollary}
\newtheorem{lem}[thm]{Lemma}

\newtheorem{prop}[thm]{Proposition}

\numberwithin{thm}{subsection} 
\endgroup

\newtheorem*{thm*}{Theorem}
\newtheorem*{conj*}{Conjecture}
\newtheorem*{verm*}{Vermutung}

\theoremstyle{definition}
\newtheorem{defn}[thm]{Definition}
\newtheorem{rem}[thm]{Remark}
\newtheorem{example}[thm]{Example}

\newtheorem{notation}[thm]{Notation}

\numberwithin{equation}{section}

\newcommand{\nosubsections}{\renewcommand{\thethm}{\thesection.\arabic{thm}}
                            \setcounter{thm}{0}
                           }

\newcommand{\cref}[3]{(\ref{#1}, #2 \ref{#3})}

\date{\today}

\newcommand{\figuretype}{eps}

\usepackage{epsfig}
\setboolean{probleme}{true}

\newcommand{\Fq}{\FZ_q}

\newcommand{\binomial}[2]{{#1 \choose #2}}

\newcommand{\secemailbothmer}{
\setlength{\unitlength}{1pt}
bothmer
\begin{picture}(0,1)
\put(0,0){m}
\put(-5,0){@}
\end{picture}
ath.uni-hannover.de}

\newcommand{\secemailschreyer}{
\setlength{\unitlength}{1pt}
schreyer
\begin{picture}(0,1)
\put(0,0){m}
\put(-5,0){@}
\end{picture}
ath.uni-sb.de}

\DeclareMathOperator{\\Bin}{\Bin}

\begin{document}

\title[A quick and dirty irreducibility Test]
      {A quick and dirty irreducibility Test for Multivariate 
       Polynomials over $\Fq$}

\author{H.-C. Graf v. Bothmer}
\address{Institut f\"ur Mathematik (C), Welfengarten 1, Universit\"at Hannover, D-30167 Hannover, Germany}
\email{\secemailbothmer}

\author{F.-O. Schreyer}
\address{Mathematik und Informatik, Geb. 27,
Universit\"at des Saarlandes,
D-66123 Saarbr\"ucken,
Germany}
\email{\secemailschreyer}

\maketitle

\begin{abstract}
We provide some statistics about an irreducibility/reducibility test
for multivariate polynomials over finite fields based
on counting points. The test works best for polynomials in a large
number of variables and can also be applied to black box polynomials.
\end{abstract}

\section{Introduction}
\nosubsections

Let $f \in \Fq[x_1,\dots,x_n]$ be a polynomial. Since $f(x)$ can take only
$q$ possible values for every
point in $x \in \AZ^n(\Fq)$ we expect that $f(x)=0$ for about
$\frac{1}{q}$ of the points $\AZ^n(\Fq)$. If on the other hand $f=gh$ is a product of
two polynomials $g,h \in \Fq[x_1,\dots,x_n]$, we have $f(x)=0$
if $g(x)=0$ or $h(x)=0$. So one might expect that products
of polynomials satisfy $f(x)=g(x)h(x)=0$ for approximately $\frac{2}{q}-\frac{1}{q^2}$
of the points $x \in \AZ^n(\Fq)$. This phenomenon is well explained by the Weil formulas \cite{etale}.

In this article we investigate the following irreducibility test for multivariate polynomials $f$
over $\Fq$:

%
Evaluate 
$f$ at $N$ random points. We reject the hypothesis that $f$ is reducible, if
the fraction of zeros $\gamma_q(f)$ found is significantly smaller than $\frac{2}{q}-\frac{1}{q^2}$.
Note that $99.5\%$ of all polynomial functions satisfy 
\[
        \gamma_q(f) \le \frac{1}{q} + 2.58 \sqrt{\frac{\frac{1}{q}(1-\frac{1}{q})}{q^n}}.
\]
This irreducibility test is quick, since the number of evaluations needed to detect
a given percentage $1-\epsilon$ of all products of polynomial functions or all general polynomial functions
do not depend on the degree of the polynomials considered respectively, i.e.
\[
   N \sim O(-q\ln \epsilon).
\]
On the other hand it is dirty, since it does not give a definite answer. 
Moreover we cannot make $\epsilon$ 
arbitrarily small, because $N$ is bounded by $q^n$, the number of $\Fq$ rational
points in $\AZ^n(\Fq)$. There will always be a few polynomials that cannot be correctly classified by our method at all. For example the product of an irreducible, absolutely reducible polynomial with
a further absolutely irreducible polynomial.

The test works for implicitly 
given (black box) polynomials as well. We give examples of such polynomials below.

The expected fraction of zeros for special classes of polynomials is also
larger than $\frac{1}{q}$. For example, the expected fraction of
zeros for $n \times n$ determinants is
\[
        E(\gamma_{q,\det}) = 1/q + 1/q^2 - 1/q^5 - 1/q^7 + O(1/q^{12})
\]
for $n \ge 12$.

\newcommand{\Bin}{\sB}
\newcommand{\Norm}{\sN}

\begin{notation} \quad

\begin{tabular}{ll}
$\Fq$               & the finite field with $q$ elements \\
$X \subset \AZ^n$   & an affine algebraic set \\
$X(\Fq)$        & the $\Fq$-rational points of $X$ \\
$|X| = |X(\Fq)|$       & the number of $\Fq$-rational points of $X$\\
$\gamma_q(X)$       & the fraction of $\Fq$-rational points in $\AZ^n$ \\
                & that
                      are contained in $X$  \\             
$\Bin(N,p,k) = {N \choose k}p^k(1-p)^{N-k}$ 
                    & the binomial distribution \\
$N$             & the number of trials \\
$p$            & the success probability\\
$k$            & the number of successes \\                    
$\Norm(\mu,\sigma)$ & the normal distribution with mean $\mu$ \\
                & and variance $\sigma^2$\\             
\end{tabular}

$\Bin(N,p)$ can be approximated by $\Norm(p,\sqrt{p(1-p)/N})$. 
\end{notation}

\section{Fractions of Zeros}
\nosubsections

\begin{example}
We choose fixed polynomials $f_1,f_2$ of degree $5$ and $f_3$ of degree
$10$ in $\ZZ[x_1,\dots,x_4]$ with coefficients in $[-9,9]$ using the random number generator of the computer
algebra system MACAULAY 2 \cite{M2} and consider $f=f_1f_2+7f_3$.
Let $X$ be the vanishing set $V(f)$.
\end{example}

A black-box polynomial is a polynomial for which  it is easy
to  check $f(x)=0$, but the explicit formula for $f$ in terms of the unknowns
$x_1\dots x_n$ is hard or impossible to write down.

\begin{example}\label{singPlaneCurves}
Let $S_d \subset H^0(\PP^2,\sO(d))$ be the hypersurface of singular homogeneous
polynomials $f$ of degree $d$ in $3$ variables. For each point $f \in H^0(\PP^2,\sO(d))$
it is easy to decide whether $f \in S_d$ via the Jacobi criterion. On the
other hand the equation of $S_d$ in the ${{d+2 \choose 2}}$ variables is not obvious.
\end{example}

\begin{example}
Let $C \subset \PP^4$ be the determinantal curve of degree $10$ and genus $6$
defined by the maximal minors
of the following $5 \times 3$ matrix
\[    
\begin{small}
\begin{pmatrix}{{x}}_{0}+{{x}}_{1}-{{x}}_{{3}}-{{x}}_{{4}}&
      {{x}}_{0}-{{x}}_{1}-{{x}}_{{2}}-{{x}}_{{4}}&
      -{{x}}_{0}+{{x}}_{{3}}+{{x}}_{{4}}\\
      -{{x}}_{0}-{{x}}_{{2}}+{{x}}_{{3}}+{{x}}_{{4}}&
      {{x}}_{0}-{{x}}_{1}-{{x}}_{{2}}-{{x}}_{{3}}+{{x}}_{{4}}&
      -{{x}}_{0}+{{x}}_{1}-{{x}}_{{2}}+{{x}}_{{3}}+{{x}}_{{4}}\\
      -{{x}}_{0}-{{x}}_{{2}}-{{x}}_{{3}}-{{x}}_{{4}}&
      -{{x}}_{0}-{{x}}_{1}-{{x}}_{{3}}-{{x}}_{{4}}&
      -{{x}}_{1}+{{x}}_{{4}}\\
      -{{x}}_{1}-{{x}}_{{2}}-{{x}}_{{3}}+{{x}}_{{4}}&
      -{{x}}_{1}-{{x}}_{{2}}&
      -{{x}}_{1}+{{x}}_{{2}}\\
      -{{x}}_{0}+{{x}}_{1}-{{x}}_{{2}}-{{x}}_{{3}}-{{x}}_{{4}}&
      -{{x}}_{0}+{{x}}_{{2}}-{{x}}_{{3}}+{{x}}_{{4}}&
      {{x}}_{0}-{{x}}_{1}+{{x}}_{{2}}+{{x}}_{{3}}+{{x}}_{{4}}\\
      \end{pmatrix}
\end{small}
\]
Let $D = \{H \in \check{\PP}^4 \suchthat H \cap C \, \text{is singular} \}$ be the dual
variety of $C$. 
\end{example}

\begin{defn}
Let $X \subset \AZ^n$ an algebraic set. We denote by
\[
       \gamma_q(X) := \frac{| X(\Fq)|}{| \AZ^n(\Fq)|}
\]
the {\sl fraction of $\Fq$-rational points} on $X$. In particular for a hypersurface
$X = V(f)$ we have $\gamma_q(f)=\gamma_q(V(f))$. We call $\gamma_q(f)$ the 
{\sl fraction of $\Fq$-rational zeros of $f$}.
\end{defn}

\begin{example}
We estimate $\gamma_q$ in three of our examples by evaluating in $N=1000$ random points over
all primes up to $17$. The following table gives the $99\%$ confidence interval for $\gamma_q$:
\[
\begin{array}{|c|r|r|r|}
\hline
q & \multicolumn{1}{|c|}{X} & \multicolumn{1}{|c|}{S_8} & \multicolumn{1}{|c|}{D}\\
\hline
2& 56.7\%\pm4.0\% & 68.4\%\pm2.9\% & 55.3\%\pm4.1\% \\
3& 33.8\%\pm3.9\% & 42.3\%\pm3.1\% & 49.2\%\pm4.1\%\\
5& 17.9\%\pm3.1\% & 24.0\%\pm2.6\% & 24.9\%\pm3.5\% \\
7& 26.2\%\pm3.6\% & 16.8\%\pm2.3\% & 35.3\%\pm3.9\%\\
11& 9.3\%\pm2.4\% &  8.9\%\pm1.8\% & 8.0\%\pm2.2\% \\
13& 8.6\%\pm2.3\% &  9.6\%\pm1.8\% & 8.4\%\pm2.3\%\\
17& 5.2\%\pm1.8\% &  8.1\%\pm1.7\% & 5.9\%\pm1.9\%\\
\hline
\end{array}
\]
In this article we will explain these numbers.
\end{example}

\begin{rem}
We can compute the true values $\gamma_2(X)=56.3\%$, $\gamma_3(X)=34.6\%$, $\gamma_5(X)=18.7\%$ and $ \gamma_7(X)=27.6\%$ with the same effort, since there are less
than $1000$ rational points in $\AZ^4(\Fq)$ for $q\le7$.
\end{rem}

\newcommand{\FqPolynomials}{\Fq[x_1\dots x_n]}

\newcommand{\functions}{R}

To study the map
\[
        \gamma_q \colon \FqPolynomials \to [0,1], \, f \mapsto \gamma_q(f)
\]
we note that $\gamma_q(f)$ factors over the ring $\functions := \map(\AZ^n(\Fq),\Fq)$:
\xycenter{
                        \FqPolynomials  \ar[r]^-{\gamma_q} \ar[d]^\psi & [0,1] \\
                        \functions \ar[ur]
                }

\begin{lem} 
$\psi$ is surjective. 
\end{lem}

\begin{proof}
Since $|\AZ^n(\Fq)| = q^n < \infty$ we can find a polynomial with prescribed values at these points
via interpolation. 
\end{proof}




We study the distribution of $\gamma_q$ on $\functions$ by regarding it as 
a random variable on the finite probability space
\[
          (\functions,\Omega,P)
\]
with $\Omega$ the sigma algebra of all subsets of $\functions$ and $P$ the constant probability measure. 

\begin{prop}
The distribution of $\gamma_q$ on $\functions$
is binomial 
\[
P\left(\gamma_q=\frac{k}{q^n}\right) = \\Bin\left(q^n,\frac{1}{q},k\right).
\]
In particular the expectation value of $\gamma_q$ is 
$E(\gamma_q) = \frac{1}{q}$. 
\end{prop}

\begin{proof}
We have to count the maps $f \in \functions$ that  map precisely $k$ different points 
to $0$. Since the values at different points are independent, this number is
\[
      \binomial{q^n}{k}
           1^k 
           \cdot 
           (q-1)^{q^n-k}
\]
The probability that $\gamma_q = \frac{k}{q^n}$ is
therefore
\[
    P\left(\gamma_q = \frac{k}{q^n}\right) = 
           \binomial{q^n}{k}
           \left(\frac{1}{q}\right)^k 
           \cdot 
           \left(\frac{q-1}{q}\right)^{q^n-k} = 
    \Bin\left(q^n,\frac{1}{q},k\right)
\]
\end{proof}

\begin{example} \xlabel{e-allF11}
Consider  maps $f \in R = \map(\AZ^4(\FF_{11}), \FF_{11})$. 
The distribution of fractions of zeros is
\[
     P\left(\gamma_{11} = k /11^4 \right) = 
     \Bin\left(11^4,1 / 11,k\right).
\]
From its approximation by the normal distribution
$\Norm(0.0909,0.0024)$ we obtain
\[
    P(0.0847 \le \gamma_{11} \le 0.0971) \le 99\%.
\]
\end{example}

We now consider products. 
The random variable
\[
       \gamma_{q,\cup} \colon \functions \times \functions
                         \to [0,1], \, \gamma_{q,\cup}(f,g) = \gamma_q(fg) = |V(f) \cup V(g))|/q^n
\]
which assigns to each pair of functions the fraction of zeros of their product.

\begin{prop}
On $\functions \times \functions$ the distribution of $\gamma_{q,\cup}$
is
\[
    P\left(\gamma_{q,\cup} = k/ q^n \right) = 
    \Bin\left(q^n,(2q-1)/q^2,k\right).
\]
In particular the expectation value of $\gamma_{q,\cup}$ is 
\[
   E(\gamma_{q,\cup}) = \frac{2q-1}{q^2}= 1 - \left( \frac{q-1}{q} \right)^2.
\]
\end{prop}

\begin{proof}
The value of $f\cdot g$ in a point $x$ depends on the values of $f$ and
$g$ at $x$. There are $q^2$ ways of choosing these values of which
$(q-1)^2$ give $(f\cdot g)(x) \not=0$.
\end{proof}

\begin{example}
Consider pairs $(f,g)$ of functions in $R$ as in Example \ref{e-allF11}.
The distribution of $\gamma_{11,\cup}$ is now
\[
     P\left(\gamma_{11,\cup} = k / 11^4 \right) = 
     \Bin\left(11^4,21/11^2,k\right).
\]
From its approximation by the normal distribution $\Norm(0.1736,0.0031)$,
we obtain
\[
            P(0.1655 \le \gamma_{11,\cup} \le 0.1816) \ge 99\%
\]
Note that this range does not intersect
\[
    P(0.0847 \le \gamma_{11} \le 0.0971) \ge 99\%.
\]

\begin{figure}
{\bf Points on a hypersurface of degree $10$ in $\AZ^4$}
\epsfig{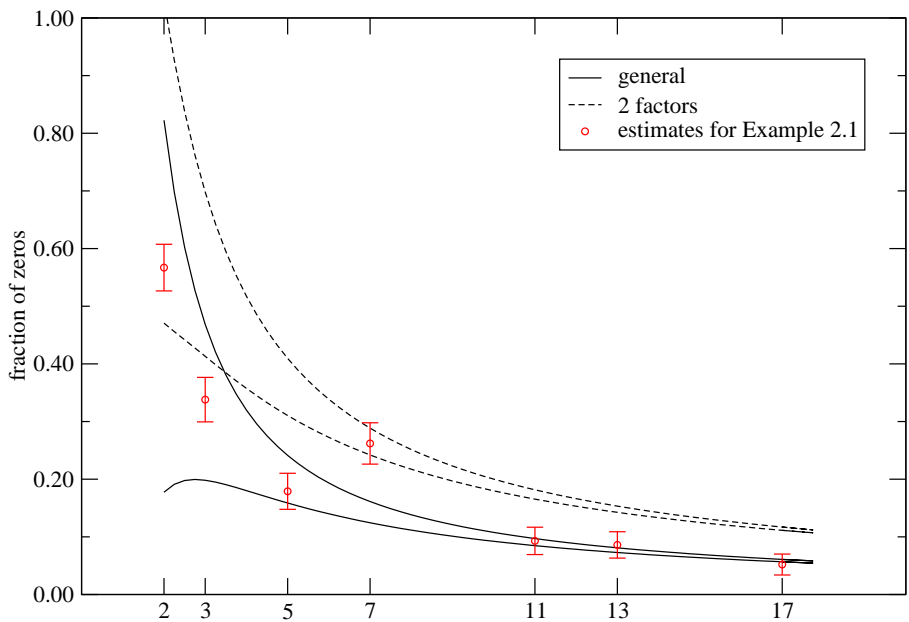}
\caption{$99\%$ of polynomial functions on $\AZ^4$ 
have $\gamma_q$ between the continuous lines. $99\%$ of products
have $\gamma_q$ between the dashed lines.}
\end{figure}

\end{example}


Geometrically products of functions correspond to the union of
their zero-sets. $\gamma_q$ also behaves well under other geometric
operations:

\begin{prop}[Intersection] \xlabel{p-intersection}
Let $X \subset \AZ^n$ be a subvariety. 
We consider the random variable
\[
        \gamma_{q,\cap X} \colon \functions  \to [0,1], \, \gamma_{q,\cap X}(f) = | V(f) \cap X|/q^n.
\]
The distribution of $\gamma_{q,\cap X}$ is
\[
        P(\gamma_{q,\cap X}=k/q^n ) = \Bin(|X|,1/q,k).
\]
In particular, the 
expectation value of $\gamma_{q,\cap X}$ is  
$E(\gamma_{q,\cap X}) = \gamma_q(X)/ q$,
where
$
        \gamma_q(X) = |X| / q^n
$
is the fraction of points on $X$ in $\AZ^n(\Fq)$.
\end{prop}

\begin{proof}
Clearly, $x \in X \cap V(f)$ if and only if $x \in X$ and
$f(x)=0$. Since the values of $f$ can be chosen independently
on the points of $X$, we have
\[
    P(x \in \ker f \cap X | x \in X) = \frac{1}{q}.
\]
\end{proof}


\begin{cor}
Consider the random variable
\[
      \gamma_{q,\cap} \colon \functions^c \to [0,1], \,
      \gamma_{q,\cap}(f_1,\dots,f_c) = | V(f_1) \cap \dots \cap V(f_c) |/ q^n
\]
Then the expected fraction of points
is
$
   E(\gamma_{q,\cap}) = \frac{1}{q^{c}}
$
\end{cor}

\begin{proof}
Use Proposition \ref{p-intersection} inductively.
\end{proof}

Notice that for polynomials $f_1,\dots, f_c$ the expected codimension
of $V(f_1,\dots,f_c) \subset \AZ^n$ is also $c$.

\begin{prop}[Substitution] \xlabel{p-substitution}
Let $\functions^m = \map(\AZ^n(\Fq),\AZ^m(\Fq))$ and $X \subset \AZ^m(\Fq)$
a subset.  Consider the random variable
\[
            \gamma_{q,subst} \colon \functions^m \to [0,1], \,
           \gamma_{q,subst}(\phi)= | \phi^{-1} X | / q^n
\]
The distribution of $\gamma_{q,subst}$ is
\[
    P\left(\gamma_{q,subst} = k / q^n \right) = 
    \Bin\left(q^n,\gamma_q(X),k\right).
\]
In particular the expectation value of $\gamma_{q,subst}$ is
$
         E(\gamma_{q,subst}) = \gamma_q(X) = |X| / q^n.
$
\end{prop}

\begin{proof}
Choosing functions $f_1,\dots,f_n$ is equivalent
to independent choice of the image points.
Therefore the probability of $\phi^{-1}(X)$ containing exactly $k$ points 
is the same
as the probability of hitting $k$ points of $X$ while choosing $q^n$ points
in $\Fq^n$. This gives the desired binomial distribution. 
\end{proof}

\section{Determinantal Varieties}
\nosubsections

Even though we have shown, that $E(\gamma_q)= \frac{1}{q}$ with
a small variance on the set of all functions from $A$ to $\Fq$, there are special classes of functions
that have larger expected $\gamma_q$. 

\begin{figure}
{\bf Singular curves in $\PP^2$}
\epsfig{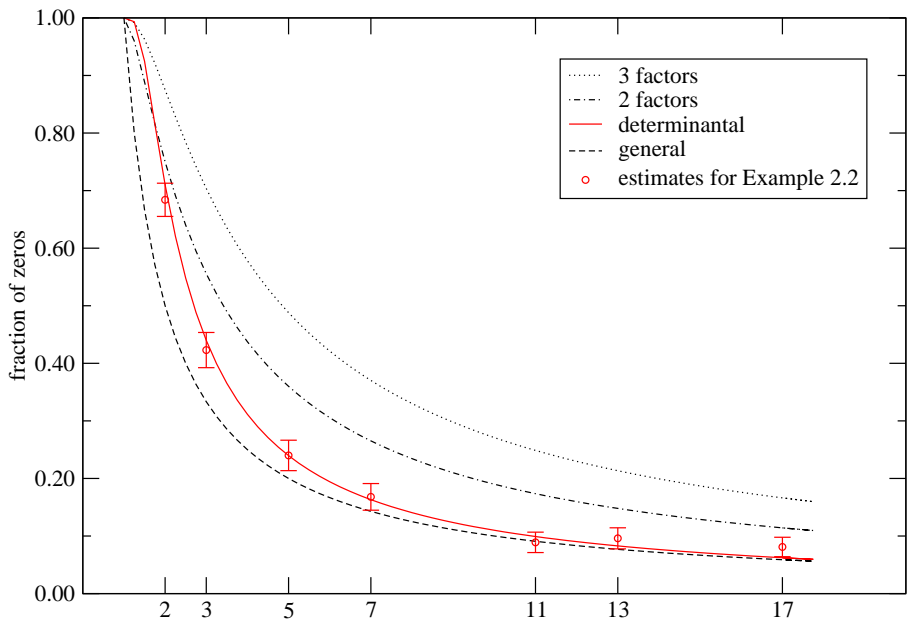}
\caption{The diagram shows the expectation values for various classes
  of polynomials in a large number of variables, and the measurement
  for $S_8$, the hypersurface of singular plane curves of degree 8.
Note that the diagram tells that about 70 \% of all plane curves over
  $\FZ_2$
are singular. 
}\label{sC}
\end{figure}

It turns out that this behavior is common for determinants:

\begin{prop} Let $X \subset \AZ^{nm}$ 
be the determinantal variety of $n \times m$ matrices with $n \le m$
of rank less than $n$.
Then the fraction of points on $X$ is
\[
    \gamma_q(X)= 1 - \prod_{i=0}^{n-1} 
                   \left( 1-\frac{1}{q^{m-i}} \right),
\]
i.e. $X$ contains $\gamma_q(X)\cdot q^{nm}$ points.
\end{prop}

\begin{proof}
We prove that the number of matrices that have maximal rank 
is
\[
     \prod_{i=0}^{n-1} \left( q^m-q^i \right)
\]
by induction.
$M$ is a matrix of full rank if and only if the first $n-1$
rows form a matrix of full rank and the last row is linearly independent
of the first $n-1$ rows. Since there are $q^{n-1}$ linear combinations
of the first $n-1$ rows we obtain a further factor $(q^m-q^{n-1})$.

\end{proof} 


\begin{cor} \label{c-det}
On the space of matrices $\functions^{nm}$,
consider the random variable
\[
         \gamma_{q,\det} \colon \functions^{nm} \to [0,1], \,
         \gamma_{q,\det}(M) = |\{ x \in \AZ^n  \suchthat \rank M(x) < n\}| / q^n.
\]
Then the fraction of zeros has expectation value
\[
    E(\gamma_{q,\det})= 1 - \prod_{i=0}^{n-1} 
                      \left( 1-\frac{1}{q^{m-i}} \right) = \frac{1}{q^{m-n+1}}+\dots
\]
The distribution of $\gamma_{q,det}$ is 
\[
    P\left(\gamma_{q,\det}= k / q^n \right) = \Bin(q^n,E(\gamma_{q,\det}),k)
\]
\end{cor}

\begin{proof}
Substitute functions for the variables in the
generic $n \times m$ matrix and use Proposition \ref{p-substitution}
\end{proof}

In the special case of $n \times n$ square matrices we have
\[
        E(\gamma_{q,\det}) = 1/q + 1/q^2 - 1/q^5 - 1/q^7 + O(1/q^{12})
\]
for $n \ge 12$.

\begin{example}[Example \ref{singPlaneCurves} continued] For small
    primes 
the divisor $S_d$ has 
more points than expected for irreducible polynomials, but not
enough to seem reducible, see Figure \ref{sC}. Our measurements are
consistent with the well known fact that $S_d$ is an irreducible 
determinantal hypersurface \cite[Chapter 13, Prop. 1.6 and 1.7]{GKZ}.
\end{example}


\begin{figure}
{\bf Points on the dual variety of a curve in $C \subset \PP^4$}
\epsfig{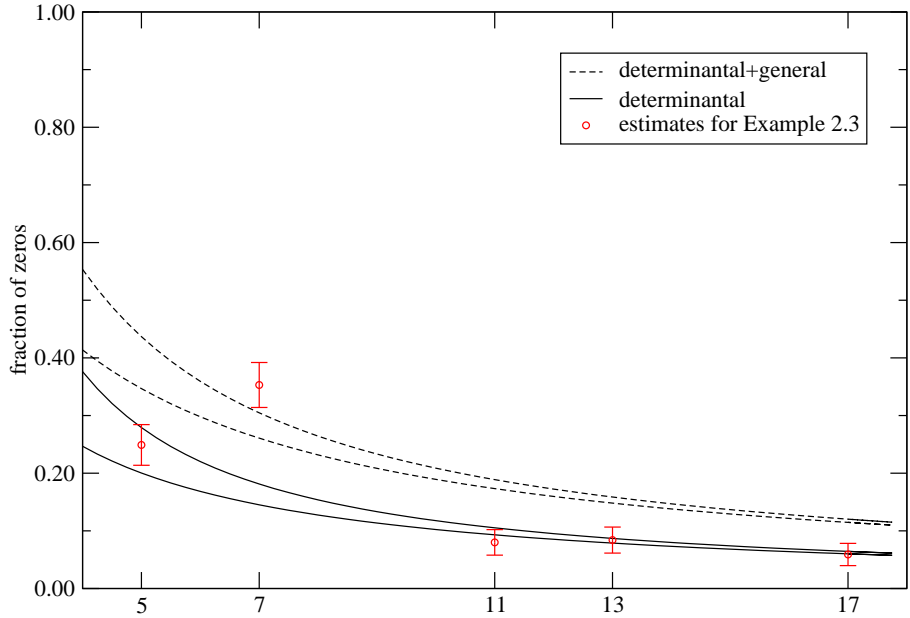}
\caption{$C$ 
has a simple node over $\FF_7$ and is smooth over $\FF_p$ for $p=5, 11, 13, 17$.}
\end{figure}






\section{Testing}
\nosubsections

To decide between two binomial distributions with success probabilities
$p_1 < p_2$ and $N$ experiments, we compute  empirical probability 
$\bar{p}=\frac{k}{N}$ and decide for $p_1$ if 
\[
    \bar{p} \le p_{middle} = \sqrt{p_1p_2} \frac{\sqrt{p_1(1-p_2)} + \sqrt{p_2(1-p_1)}}
    {\sqrt{p_1(1-p_1)} + \sqrt{p_2(1-p_2)}} \approx  \sqrt{p_1p_2}.
\]
To achieve a confidence level of $1-\epsilon$
we choose $s = s(\epsilon)$ such that 
\[
        \Phi(s) = \frac{1}{\sqrt{2\pi}} \int_s^\infty e^{-\frac{x^2}{2}}dx = \epsilon
\]
and $N$ such that 
\[
       \sqrt{N} \ge s(\epsilon) \frac{\sqrt{p_1(1-p_1)}+\sqrt{p_2(1-p_2)}}{p_2-p_1}.
\]
In our case we have
\[
     p_1 \le \frac{1}{q} + s(\epsilon) \sqrt{\frac{\frac{1}{q}(1-\frac{1}{q})}{q^n}}
\]
for $1-\epsilon$ of all polynomials and
\[
     p_2 \ge \frac{2q-1}{q^2}- s(\epsilon) \sqrt{\frac{\frac{2q-1}{q^2}(1-\frac{2q-1}{q^2})}{q^n}}.
\]
for $1-\epsilon$ of all products of polynomials. The decision based on the empirical 
probability $\bar{p}=\frac{k}{N}$, is then correct in $1-\epsilon$ cases of the experiments.
Note however, that
for fixed $n$ and $q$ we cannot make $\epsilon$ arbitrarily small, since we need $p_1 \le p_2$.

An easy calculation gives the following estimate
\[
        \sqrt{N} \ge s(\epsilon) \frac{(2q)^\frac{3}{2}}{q - 1- 2s q^{-\frac{n-2}{2}}}
\]
for $q \ge 3$, which approaches $s(\epsilon)\sqrt{2q}$ for large $n$ or $q$. 
Since $s(\epsilon) = O(\sqrt{-\ln(\epsilon)})$, we conclude that $N$ grows like
$O(-q\ln \epsilon)$.

For $\epsilon = 0.5\%$, $s=2.58$, the number of trials needed is
\[ 
 \begin{array}{|c|rrrrrrr|} 
\hline
       & 2 & 3 & 5 & 7 & 11 & 13 & 17 \\
\hline  
   n=1 & \infty &\infty &\infty &\infty &\infty &\infty &\infty \\
   n=2 & \infty &\infty &\infty & \infty &\infty &\infty &\infty \\
   n=3 & \infty &\infty & \infty &28373&2355&1908&1669\\
   n=4 & \infty &\infty &1103&647&634&682&803\\
   n=5 & \infty &1705&367&369&482&551&695\\
   n=6 & \infty &384&259&308&447&521&673\\
   n=7 & 4457&224&225&289&437&513&667\\
   n=8 & 619&173&212&283&434&511&666\\
   n=9 & 295&151&206&280&433&511&666\\
   n=10 & 197&140&204&279&433&511&665 \\
\hline
\end{array}
\]
$\infty$ indicates that there are not enough points in $\AZ^n(\Fq)$ to
perform the test for the required $\epsilon= 0.5 \%$.
In  case we can perform the test,  the deciding number of successes $Np_{middle}$ is
\[ 
\begin {array}{|c|rrrrrrr|} 
\hline
       & 2 & 3 & 5 & 7 & 11 & 13 & 17 \\
\hline  
   n=1 & \infty &\infty &\infty &\infty & \infty &\infty &\infty \\
   n=2 & \infty &\infty &\infty &\infty &\infty &\infty &\infty \\
   n=3 & \infty &\infty &\infty &5607&301&207&139\\
   n=4 & \infty &\infty &303&128&81&74&66\\
   n=5 & \infty &754&101&73&61&59&57\\
   n=6 & \infty &170&71&61&57&56&55\\
   n=7 & 2821&99&61&57&55&55&55\\
   n=8 & 391&76&58&56&55&55&55\\
   n=9 & 186&67&56&55&55&55&55\\
   n=10 & 125&62&56&55&55&55&55\\
\hline
\end {array}
\]

\section{Higher codimension}
\nosubsections

\begin{figure}[t] 
{\bf Surfaces in  $\PP^4$}
\epsfig{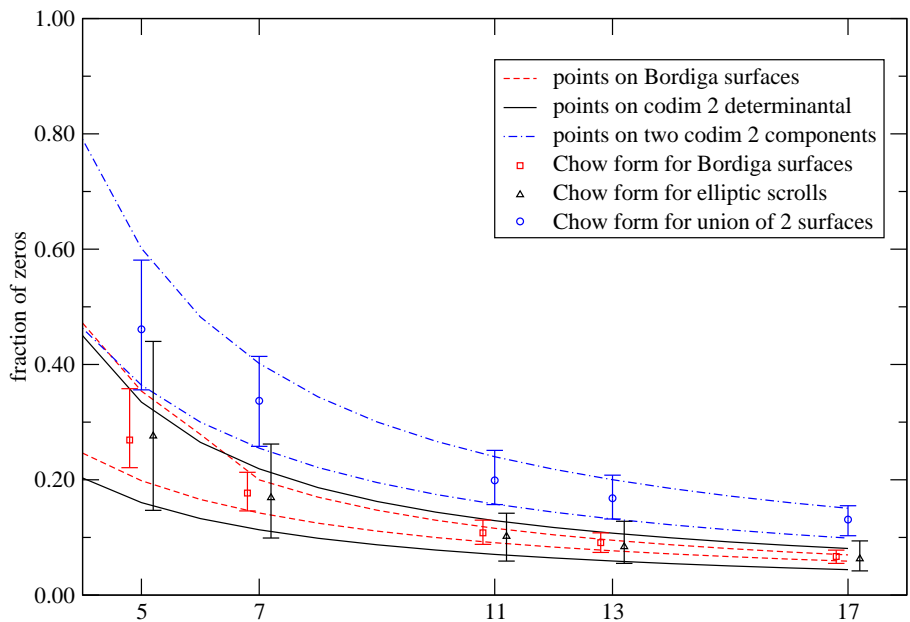}
\caption{The $5\%$ and
the $95\%$ quantiles of  $\gamma_q$ for the Chow forms of 
$100$ Bordiga surfaces, elliptic scrolls and their unions compared
with the error estimates for counting points on codimension $2$ determinantal
varieties rescaled. Using the geometry of Bordiga surfaces we obtain a better estimate.
}\label{f-boot}
\end{figure}

In principle this method can be applied to algebraic sets of higher codimension. 

Consider two surfaces in $\PP^4$ and their union. We would like to distinguish their union form
the irreducible examples. One possibility is to consider the Chow form which is a determinantal
hypersurface on $G(2,5)$ in this case. In Figure \ref{f-boot} we indicate the $5\%$ and
the $95\%$ quantiles of  $\gamma_q$ for the Chow forms of 
$100$ Bordiga surfaces, elliptic scrolls and their unions.
 A second possibility is to  count points and apply Corollary \ref{c-det}. 
As Figure \ref{f-boot} shows there is 
no difference between the two methods. The formula for the error term
underestimates 
the number of points on a elliptic scroll,
because the scroll is irregular.

The method of searching points at random in higher codimensional subsets of
rational varieties helped us in proving the existence of several
interesting components of Hilbert schemes. \cite{smallFields}, \cite{needles}, \cite{newfamily}

\def\cprime{$'$}


\end{document}